\documentclass[12pt]{article}

 \usepackage{hyperref}
\usepackage{graphicx}
\usepackage{amssymb}
\usepackage{epstopdf}
\usepackage{amsmath}

\newcommand{\ggamma}{\boldsymbol{\gamma}}

\newcommand{\flambda}{\frac{1}{\lambda}}
\newcommand{\ex}{e^{i\theta}}
\newcommand{\mex}{e^{-i\theta}}
\newcommand{\zbar}{\bar{z}}

\newcommand{\parbar}{\bar{\partial}}
\newcommand{\fracpar}[2]{\frac{\partial #1}{\partial #2}}
\newcommand{\be}{\begin{equation}}
\newcommand{\ee}{\end{equation}}
\newcommand{\bea}{\begin{eqnarray}}
\newcommand{\eea}{\end{eqnarray}}
\newcommand{\eps}{\epsilon}
\newcommand{\C}{\mathbb{C}}
\newcommand{\R}{\mathbb{R}}

\newcommand{\x}{\mathbf{x}}

\newtheorem{theorem}{Theorem}[section]
\newtheorem{lemma}[theorem]{Lemma}

\newtheorem{proposition}[theorem]{Proposition}
\newtheorem{corollary}[theorem]{Corollary}

\newenvironment{proof}[1][Proof]{\begin{trivlist}
\item[\hskip \labelsep {\bfseries #1}]}{\end{trivlist}}
\newenvironment{definition}[1][Definition]{\begin{trivlist}
\item[\hskip \labelsep {\bfseries #1}]}{\end{trivlist}}

\newenvironment{remark}[1][Remark]{\begin{trivlist}
\item[\hskip \labelsep {\bfseries #1}]}{\end{trivlist}}

\newcommand{\qed}{\nobreak \ifvmode \relax \else
      \ifdim\lastskip<1.5em \hskip-\lastskip
      \hskip1.5em plus0em minus0.5em \fi \nobreak
      \vrule height0.75em width0.5em depth0.25em\fi}

\title{Ray Transforms and Vector Fields}
\date{February 7, 2011}
\author{Nicholas Hoell}
\begin{document}
\maketitle

\begin{abstract}
\noindent We review and extend a technique for recovering a smooth function from its averages over a wide class of curves in a general region of Euclidean space.  The method is based on complexification of the underlying vector fields defining the transport and recasting the problem in terms of complex-analytic function theory.  Conditions on the validity of prior formulae appearing in \cite{hoell} as well as stability estimates are then discussed first for the case of vector fields with polynomial coefficients and later for more general cases. 
\end{abstract}

\section{Background and Motivation}
The following filtered backprojection formula appeared in \cite{Bal1};
\be
\label{bals}
f(z)=\frac{(1-|z|^2)^2}{4\pi}\int_0^{2\pi}\frac{1}{|1-z\mex|^4} H_s \fracpar{}{s}  I f(s(z\mex),\ex) d\theta  
\ee
where $f(z)$ is a smooth enough function and $If$ its geodesic ray transform in the Poincar\'e disc.  Formula (\ref{bals}), once obtained, subsequently gave a \textbf{holomorphic integrating factor} to derive a similar, though more complicated, formula for the \textbf{attenuated radon transform} (AtRT) $I_a f$ on the same space.  The method used in that paper, which we call \textit{the method of complexification}, was an extension of one first used in \cite{novikov} and rests on the introduction of a complex parameter $\lambda$ into the governing transport equation and a subsequent analysis of the behavior of solutions in terms of this new parameter.  

 Recently, in \cite{hoell}, we obtained a strikingly similar result to (\ref{bals}), namely 
\be
\label{mine}
f(z)=\frac{1}{4\pi} \int_0^{2\pi}P(\lambda_i,\theta)X^\bot_\theta H(I_\theta f)(s(z\mex),\ex)d\theta
\ee
The aim of this paper is to outline, briefly, the method which resulted in the above formula as well as to further classify the breadth of its validity.  Although the original impetus for the above was the inversion of the AtRT, formula (\ref{mine}) is interesting in its own right.  Much of the material in this article may be found in more detail in \cite{hoell2}

The structure of this paper is as follows.  Section \ref{preliminaries} contains the entire cast of characters needed throughout the article.   Since our goal is partly expository, section \ref{complexification} is a concise review of the method of complexification as well as a discussion of what we term ``H-ness" and its limitations. In section \ref{generalities} we examine H-ness in more detail first for the case of polynomials and eventually for more general real-analytic vector fields.  We generalize the situation in the penultimate section \ref{harmonic} where we present results relating H-ness to reasonable frequency constraints.   

\section{Generic Preliminaries}\label{preliminaries}
Our setup will be as follows.  Let $\ggamma: \R^2 \ni (t,s) \mapsto \ggamma(t,s) \in \Sigma \subset \R^2$ be a real-analytic diffeomorphism on a simply-connected domain $\Sigma$ generating the linear, stationary transport operator
\be
 X_{\x} =a(x,y)\fracpar{}{x}+b(x,y)\fracpar{}{y}, \qquad \x=(x,y)\in \Sigma \ee
We regard $\R^2 \cong \C$ via the standard isomorphism so that $\ggamma$ is identified with $\gamma^1(t,s)+i\gamma^2(t,s)$. Defining complex $w\doteq \ggamma(t,s)$, we see that $(w,\bar{w})$ are now (independent) \textit{complex} coordinates on $\Sigma$.   The regularity of the curves $\ggamma(t,s)$ show us that $\ggamma_* \fracpar{}{t}$ is a non-degenerate field on $\Sigma$, $X|_w=\mu(w)\fracpar{}{w} +\bar{\mu}(w)\fracpar{}{\bar{w}}$ where $(\phi_* X)(f)=X(\phi^* f)$ defines the pushforward $\phi_*$.  

The equation of interest is the stationary transport boundary value problem $X|_w u(w)=f(w)$, for $w \in \Sigma$, $f(w) \in C^\infty _0(\Sigma)$ with $\lim_{t\searrow -\infty} u(w(t,s))=0$, i.e. the BVP
\begin{align}
\mu(w)\fracpar{u}{w}+\bar{\mu}(w)\fracpar{u}{\bar{w}}&=f(w), \quad w \in \Sigma\\
u\left.\right|_{\partial_- \Sigma}&=0
\end{align}
One key difference between this formulation of the problem and that considered in \cite{Bal1} is that there is, a priori, no immediately obvious object to ``complexify" since we no longer have a parameterization of the points of impact on $\partial \Sigma$ as was the case then.  To circumvent the aforementioned difficulty, we appeal to the Riemann mapping theorem (\cite{nehari, krantz}) which guarantees a unique biholomorphism $z:\Sigma \to D^+$ satisfying $z(\zeta)=0$, $z'(\zeta)>0$ for $\zeta \in D^+$, where $D^+$ is the unit disc $\{z \in \C; |z|<1\}$.  Since the Riemann mapping is conformal it is necessarily (infinitesimally) factorable (as in e.g. \cite{pal}) with respect to the subvarieties of integral curves of $ X_z$.  Because of this equivalence between our initial domain $\Sigma$ and the unit disc all further results will henceforth be presented in the disc.

Since $\ggamma^* z$ maps $\R^2$ into $D^+$ we use $(z,\bar{z})$ as our coordinates on $D^+$ and have a new vector field $X|_z =z_*X|_{z(w)}$ where $\mu \mapsto \{z_*\mu\} \fracpar{z}{w} \circ z^{-1}$ and likewise for $\bar{\mu}$. By a forgivable abuse of notation we denote $\{z_*\mu\} \fracpar{z}{w} \circ z^{-1}$ by $\mu(z)$ and $\{z_*\bar{\mu}\} \fracpar{\bar{z}}{\bar{w}} \circ z^{-1}$ by $\bar{\mu}(z)$ so that 
\[\left. X\right|_z=\mu(z)\fracpar{}{z} +\bar{\mu}(z) \fracpar{}{\bar{z}}  ,\qquad z\in D^+, \qquad |\mu|>0\]
is our governing differential operator. Defining $t(z) =z_*w_*t$ and $s(z) =z_*w_*s$, smooth functions on $D^+$, the method of characteristics gives the solution to the BVP $X|_zu(z)=f(z)$, $u(z(-\infty,s))=0$  as 
\be
u(z)=(D_1 f)(z)\doteq \frac{1}{2}\int_\R f(z(t_0,s))sign(t(z)-t_0)dt_0
\ee
and since $\ggamma^* z^*=(z \circ \ggamma)^*$. we define the X-ray transform of a function $f(z)$ over the integral curves of $X|_z$, indexed by the transverse parameter $s$, to be
\be
(If)(s)=\int_\R f(z(t,s))dt
\ee

The main players we need at our disposal are as follows;
\begin{description}
\item[Symmetric Beam Transform]\footnote{N.B. $\ex$ and $\theta$ will be used interchangeably, their meaning clear from context}  \[(D_\theta \psi )(z)\doteq \frac{1}{2}\int_\R \psi (\ex z(t_0,s(z \mex )))sign(t(z\mex)-t_0)dt_0, \qquad \psi \in L^1(D^+) \]
\item[Ray Transform] \[(I\psi)(s,\ex) = (I_\theta \psi)(s)\doteq \int_\R \psi (\ex z(t,s))dt, \qquad \psi \in L^1(D^+)\]
\item[Hilbert Transform] \[(H \psi)(x) \doteq \frac{1}{\pi} p.v. \int_\R \frac{\psi(y)}{x-y}dy, \qquad \psi \in L^p(\R), \qquad p>1\]
\end{description}
as well as the classical Poisson kernel $P(z,\theta)=\frac{1-|z|^2}{|1-\mex z|^2}$ which for $z \in D^+$, $\theta \in T$ generates the harmonic functions on the unit disc.  Occasionally we will use the nonstandard notation $I^X_\theta f$ to indicate the line integral of $f$ over the integral curves of the vector field $X$.  

\section{Complexification in a Nutshell}\label{complexification}
The main result of \cite{hoell} is the following.   

\begin{theorem} \label{mainy}Under suitable conditions on $\left. X \right|_z$ and $s(z)$ there exists a function $\lambda_i(z)$ on $D^+$ such that 
\be
\label{recon} 
f(z)=\frac{1}{4\pi} \int_0^{2\pi}P(\lambda_i,\theta)X^\bot_\theta H(I_\theta f)(s(z\mex),\ex)d\theta
\ee
provides a reconstruction for the function $f$ based on the data $I_\theta f$ of ray transforms of $f$ over the integral curves of $X_\theta=\theta_* (\left. X \right|_z)$.
\end{theorem}
We now review the method outlined in that paper which was used to obtain this result.

%\noindent The heuristic schema through which the above was reached consists of the following:

%\begin{enumerate}
%\item \textbf{Model:}  Writing down the transport BVP for the dynamics
%\item \textbf{Symmetrizing:} Introducing a rotation parameter $\lambda=\ex$ into the integral curves of the model PDE
%\item \textbf{Symmetry-Breaking:} Complexifiying the parameter introduced in step 2 by moving $\lambda$ ``off-shell", i.e. $|\lambda| \neq 1$
%\item \textbf{Analysis and Asymptotics:} Evaluating the dependence of solutions to the new \textit{complexified} equation on our parameter $\lambda$ and examining limiting behavior
%\item \textbf{Reconstruction:} Using holomorphicity of the solutions to write the inversion formulae as Poisson integrals of their asymptotic boundary values found in step 4
%\end{enumerate}
%The reader may benefit from keeping the above in mind throughout the subsequent exposition.  The purpose of this article is to expand on the ``sufficient conditions" mentioned in Theorem \ref{mainy}.  

%\section{Complexification and Condition H}
\subsection{Symmetrizing and Symmetry-Breaking}
Let $\lambda = \theta \in T\doteq \partial D^+$ and define the conformal map $\lambda:(z,\bar{z}) \to (\lambda z, \flambda \bar{z})$.  If, for each $s$, $\Phi(\cdot,s)$ is a set of integral curves of $D^+$, then $z^{-1}(\lambda^* \Phi(\cdot,s))$ are conformally related curves in $\Sigma$. For $\lambda \in D^+/\{0\}$ we then consider $X_\lambda \doteq \lambda_*X|_z$ to be the so-called \textit{``complexification of $\left. X \right|_z$"}. Explicitly, $\lambda_* X|_z$ takes the form $\mu(\frac{z}{\lambda},\lambda \bar{z})\lambda \fracpar{}{z}+\bar{\mu}(\frac{z}{\lambda},\lambda \bar{z}) \flambda \fracpar{}{\bar{z}}$ or $X_\lambda =  \xi(z,\lambda) \fracpar{}{z} +\rho(z,\lambda) \fracpar{}{\bar{z}}$ with $\flambda \xi(z,\lambda)= \mu(z,\lambda) \doteq \lambda_*\mu(z)$ and $\lambda \rho(z,\lambda)=\bar{\mu}(z,\lambda)=\lambda_*\bar{\mu}(z)$.  

Define $X_\lambda^{\bot}=\pm i( -\xi(z,\lambda) \fracpar{}{z} +\rho(z,\lambda) \fracpar{}{\zbar})$ as a vector field orthogonal to $X_\lambda$ when $\lambda =\ex$.  Namely, $X_\theta \cdot X^\bot_\theta=\pm i(|\xi(z,\ex)|^2-|\rho(z,\ex)|^2)=0$ in the standard inner product $\cdot :\C^2 \to \C$. The prefactor of $i$ makes $X^\bot_\theta u(z,\theta)$ real-valued and the choice of $\pm$ is determined by whichever satisfies the condition $X_1^\bot s>0$. Because $X_1^\bot=a(z)z_*\fracpar{}{s}$ for some real-valued $a(z)$, this determines $X^\bot_1$ uniquely and since we could just as well reparameterize with $-s$ we will, without losing generality, avoid keeping track of signs by assuming $X_\lambda^\bot=i(-\xi(z,\lambda)\fracpar{}{z}+\rho(z,\lambda)\fracpar{}{\zbar})$. Likewise, define $s(z,\lambda)$ and $t(z,\lambda)$ respectively as $\lambda_*s(z)$ and $\lambda_*t(z)$ for $\lambda \in D^+/\{0 \}$.  \footnote{For functions $k(z,\lambda)$, $\fracpar{k}{z}$ and $k_z$ are equivalent, as are $\fracpar{k}{\zbar}$ and $k_{\zbar}$, and we will use them interchangably.}

\subsection{Analysis, Asymptotics, and H-ness} Our complexified transport equation now reads as follows;
\be
\label{trans}
 X_\lambda u(z,\lambda)=f(z), \qquad \lambda \in D^+
 \ee
where it should be stressed that the parameter $\lambda$ has no obvious relation to the original particle transport that started this rigmarole.  The method used to obtain Theorem \ref{mainy} involves solving equation (\ref{trans}) for $u(z,\lambda)$ and showing analytic dependence of the solution on this parameter, i.e. $\partial_{\bar{\lambda}}u(z,\lambda)=0$.  

A restricted class of vector fields known as \textbf{type H} was identified which ensure the ensuing steps work out as we need.  The following is a revised version of that definition better suited to the purposes of this paper. \footnote{Notice that the maximum principle causes the third condition to follow automatically if, in the second condition, $D^+/\{0\}$ is replaced by $D^+$ which was the situation examined in \cite{hoell}.  Since the first and last of the above are in agreement with those already considered, we need only check that the middle conditions square with what we need.}  
\begin{definition}\label{condition}A real vector field $X|_z$, complexified in the manner above \[\lambda_*X|_z=a(z,\lambda)\fracpar{}{z}+b(z,\lambda)\fracpar{}{\zbar}, \quad \lambda\in D^+/\{0\}\] is said to be of \textbf{type H} if the following holds: \begin{enumerate}
\item $a(z,\lambda)$ is a holomorphic function of $\lambda$ for $\lambda \in D^+$ and has at least one zero  $\lambda=\lambda_i(z)\in D^+$
\item $b(z,\lambda)$ is a meromorphic function of $\lambda$ for $\lambda \in D^+$ and has no zeroes in $D^+/\{0\}$
\item $\frac{a(z,\lambda)}{b(z,\lambda)}$ is a holomorphic function of  $\lambda \in D^+$ 
\item $s(z,\lambda)$, $\fracpar{s(z,\lambda)}{z}$, $\fracpar{s(z,\lambda)}{\zbar}$ are meromorphic functions of $\lambda$ for $\lambda \in D^+$
\end{enumerate}
where, as in the above, $s(z,\lambda)=\lambda_*s(z)$ is the complexified parameter specifying transverse foliation of the integral curves of $X_\lambda$.

\end{definition}
This \textbf{condition H} is the ``suitable condition" mentioned in \ref{mainy} and we will assume our vector field is of this type (i.e. $a=\mu$ and $b=\rho$).   The $\lambda_i(z)$ mentioned previously are the zeroes of the complexified $\fracpar{}{z}$ coefficient of our initial field.  Note that \textbf{condition H} is strong insofar as \textit{holomorphy itself} is a rather stringent condition.   The above criteria will heretofore be called ``H-ness".  
\subsection{A Proof Sketch}
We give a scandalously brief sketch of the proof leading to (\ref{recon}), highlighting where H-ness comes into play.  First of all, by the third condition in \ref{condition} of nondegeneracy we see that the Jacobian $\partial s(z) =|s_z(z,\lambda)|^2-|s_{\zbar}(z,\lambda)|^2 \neq 0$ holds on $\lambda \in D^+/\{0\}$ since the inequality
\[0\neq |\fracpar{s(z,\lambda)}{z}\fracpar{t(z,\lambda)}{\zbar}-\fracpar{s(z,\lambda)}{\zbar}\fracpar{t(z,\lambda)}{z}| \leq |s_z(z,\lambda)|( |t_{\zbar}(z,\lambda)|+|\frac{\xi(z,\lambda)}{\rho(z,\lambda)} || t_z (z,\lambda)|) \]
guarantees that $|s_z(z,\lambda)|^2 \neq 0$ on that same region.

We may therefore make a change of variables in $s$ to get $s_*X_\lambda=s_*X_\lambda \bar{s}(z,\lambda) \fracpar{}{\bar{s}}$, whereby our  fundamental equation $X_\lambda G_\lambda (z;z_0)=\delta(z-z_0)$ is solved explicitly by 

\begin{equation}
\label{green}
G_\lambda(z,z_0)=-\lambda \frac{ \left.  \frac{\partial(t,s)}{\partial(z,\zbar)} \right |_{z_0} }{\pi (s(z)-s(z_0))} \qquad \lambda \in D^+/\{0\} 
\end{equation}

Checking against a bump function extends this to hold weakly at $\lambda \to 0$ and density shows that $u(z,\lambda)$ is holomorphic in $\lambda$ as needed.  A similar argument works on $u_z (z,\lambda)$ and $u_{\zbar}(z,\lambda)$ by invoking the final condition of H-ness.  
From here, Hilbert's relations on the boundary values of complex-analytic functions become viable as the following result shows.
\begin{proposition}\label{prop}{\[u_{\pm}(z,\ex) \doteq \lim_{D^\pm \ni \lambda \to \ex}u(z,\lambda) = \mp \frac{1}{2i}(HI_\theta f)(s(e^{-i\theta}z),\theta)+(D_\theta f)(z)\]where the Hilbert transform $H$ is taken with respect to the first variable.} 
\end{proposition}
The proof of this comes from an explicit analysis of \ref{green} with $\lambda = 1-\eps$ and deriving the relation
\begin{equation}
\label{mess}
X_1 i s'(z,1)=\frac{1}{2} \frac{  (\left.\fracpar{}{\lambda} \frac{\xi}{\rho}) \right |_{\lambda =1} }{\frac{\xi(z,1)}{\rho(z,1)}}X^\bot_1 s(z,1) 
\end{equation}
%At this point, the fact that $\frac{1}{2} \frac{  (\left\fracpar{}{\lambda} \frac{\xi}{\rho}) \right |_{\lambda =1} }{\frac{\xi(z,1)}{\rho(z,1)}}>0$ does \textit{not} follow immediately since $\frac{\xi}{\rho}$ is no longer necessarily a Blaschke product, but may rather be a Weirstrass one.

It can be shown that by invoking the third condition of \ref{condition} we have
\[sign(is'(z,1)-is'(z_0,1))=sign(t(z,1)-t(z_0,1))\] 
from which the Sokhotskyi-Plemelj formula allows us to obtain the advertised proposition.  

\subsection{Reconstruction} We use Proposition \ref{prop} together with the classical representation of complex-analytic functions on the unit disc.  By definition of $\lambda_i$ and $X^\bot_\theta$, one has $i X_{\lambda_i}u(z,\lambda_i)=X^\bot_{\lambda_i}u(z,\lambda_i)$ so that on equating real and imaginary parts we have
 \be
 \frac{1}{2\pi} \int_0^{2\pi}P(\lambda_i,\theta)X^\bot_\theta (D_\theta f)(z)d\theta=-\frac{1}{4\pi} \int_0^{2\pi}P(\lambda_i,\theta)X_\theta H(I_\theta f)(s(z\mex),\ex)d\theta \nonumber
 \ee
 and
  \be
 f(z)=\frac{1}{4\pi} \int_0^{2\pi}P(\lambda_i,\theta)X^\bot_\theta H(I_\theta f)(s(z\mex),\ex)d\theta
 \ee
which is the result we sought.  Notice that $H$ always denotes the Hilbert transform with respect to the $s$ variable.  

\section{Statement of Results}
The goal of this paper is to establish, over the next two sections, the following result (viz. Theorem \ref{last} and Corollary \ref{main2} respectively).  

\begin{theorem}\label{maintheorem} Let $\eps >0$ be given.  Suppose that $\mu(z,\zbar)=\sum_{p+q\geq 0}a_{pq} z^p \zbar^q$ is a real-analytic function on $D^+$ and that $f\in C_c^\infty(D^+)$. Define $c_r(z,\zbar)\doteq \sum_{q-p=r}a_{pq}z^p\zbar^q$.  Let $l(z)$ and $k(z)$ be the max and min respectively of the $j$ such that $c_j(z)\neq 0$.  Suppose that $\lambda_* s$, $\lambda_* s_z$, and $\lambda_* s_{\zbar}$ are meromorphic for $\lambda \in D^+$.   
\begin{itemize}
\item If there are only finitely many nonzero $c_j(z)$, and both $l(z)+k(z)+2\geq 0$ and $0<|c_k(z)|<|c_l(z)|$ holds for all nonzero $z \in D^+$, then  there exists a vector field $Y^\eps_\theta$ such that we have a perfect reconstruction
 \[f(z)=\frac{1}{4\pi} \int_0^{2\pi}P(\lambda^\eps_i,\theta)Y^{\eps \bot}_\theta H( I_\theta f)(s(z\mex),\ex)d\theta\]

\item If, for all $z \in D^+/\{0\}$, there exist infinitely many $j \in \mathbb{Z}$ such that $c_j(z)$ and $c_{j+1}(z)$ are both nonzero, and if $\limsup_{j\to \infty} |\frac{c_{j+1}(z)}{c_j(z)}|<1$ then there exists a vector field $Y^\eps_\theta$ such that 
\[||I^{Y^\eps}f-I^Xf||_{L^q(D^+)}<C\eps \qquad 1\leq q\leq \infty \]
If, in addition, we have the Fr\'echet bound $||H_{\tilde{s}}I_\theta^{\tilde{X}}f-HI_\theta f||_{\mathcal{S}}<\delta(\eps)$ then there exists a function $\lambda_i:z\to \lambda_i(z)$ satisfying
\be
\label{ineq1}
 \sup_{z\in D^+}\left| f(z)-\frac{1}{4\pi} \int_0^{2\pi}P(\lambda^\eps_i,\theta)Y^{\eps \bot}_\theta H( I_\theta f)(s(z\mex),\ex)d\theta \right|  \leq C(\delta(\eps)) \ee where $HI_\theta f$ is the Hilbert transform in the $s$ variable of the trace of $f$ over the integral curves of $\theta_* (\mu \partial +\bar{\mu} \parbar)$

\end{itemize}
\end{theorem}

Clearly (\ref{ineq}) generalizes (\ref{recon}) in the sense that it allows an \textit{approximate} reconstruction, to arbitrary accuracy, for a large class of vector fields.  We will, in the sequel, prove Theorem \ref{maintheorem} and help explain just how broad its applicability is.  The case of polynomial fields is as good as one could hope for.  Our methodology is to establish results first for the case of polynomial vector field coefficients and later to reinterpret the terms $c_j(z)$ as the frequencies of the complexified vector field's coefficients in the case of non-polynomial fields.

\section{Polynomial Vector Fields}\label{generalities}
 \subsection{The polynomial space $\Gamma(\Omega)$}Consider the case in which $\mu(z)$ is a nonvanishing polynomial, i.e.  $\mu(z,\zbar)=\sum_{p+q\geq 0}^N a_{pq}z^p\zbar^q$ for $z \in \Omega \supset \{0\}$.
The \textit{complexified} coefficients of the field $X_\lambda$ around $\lambda=0$ are then
\[\xi(z,\lambda)=\sum_{p+q \geq 0} b_{pq}(z)\lambda^{q-p+1},\qquad \rho(z,\lambda)=\sum_{p+q\geq 0} d_{pq}(z)\lambda^{p-q-1}\]
with $ b_{pq}=a_{pq}z^p\zbar^q$ and $d_{pq}=\bar{b}_{pq}$. In order for H-ness to hold\footnote{For the moment we will be ignoring any possible problems with $s(z,\lambda)$} we will need that $\partial_{\bar{\lambda}}\xi(z,\lambda)=0$, which a priori we do not have since $q-p+1$ may very well be negative.  If $q-p+1\geq 0$ for all $(p,q)$-pairs then we are (provided we have roots and the rest of \textbf{condition H}) in the position of the previous section.  If not, i.e. \textit{if $p>q+1$ holds for some $(p,q)$-pair}, then we proceed as follows.

First of all,  we will mostly be using the \textit{local irreducible exponents} $k$, $l$ given by 
\be
\label{exponents}
k(z)\doteq \min_{c_j(z)\neq 0}(j) \qquad \text{and} \qquad l(z)=\max_{c_j(z) \neq 0}(j)
\ee
where \[c_r(z,\zbar)\doteq \sum_{\substack{p,q \\q-p=r}}a_{pq}z^p\zbar^q \]
so that the Laurent expansion of $\mu(z,\lambda)$ around $\lambda =0$ is given by $\sum_{r=k}^{l} c_r(z,\zbar)\lambda^r$.  Occasionaly we will need the \textit{global exponents} defined as
\be 
\label{k}
k_{\mu}\doteq k(\mu)=\min_{z \in \Omega}k(z) \qquad \text{and} \qquad l_{\mu}\doteq l(\mu)=\max_{z\in \Omega}l(z)
\ee
%\begin{remark}A quick calculation shows that
%\be
%c_r(z,\zbar)=\left\{ \begin{array}{1 1} z^{|r|}(\kappa(r)+\sum_{j>0}\kappa_j (r)|z|^{2j}) &\quad r\leq 0 \\
%\zbar^{r}(\eta(r)+\sum_{j>0}\eta_j (r)|z|^{2j}) &\quad r> 0 \end{array} \right
%\ee
%for some constants $\kappa(r)$, $\eta(r)$, $\kappa_j(r)$, $\eta_j(r)$.  We may abbreviate the expressions $(\kappa(r)+\sum_{j>0}\kappa_j (r)|z|^{2j})$ by $g_r(|z|^2)$.  Thus, the algebraic variety $\{z;  c_r(z,\zbar)=0\}$ defines a (possibly degenerate) circle.  
%\end{remark} 

 Obviously $-N\leq k(\mu)$ and similarly $0 \leq l(\mu) \leq N$\footnote{If $l(\mu)<0$ then we should use the complementary complexification $\lambda:(z,\bar{z}) \to (\flambda z, \lambda \bar{z})$, for $\lambda \in D^+$ and get a holomorphic $\flambda \lambda_* \mu$.  Since this is a situation which was dealt with in the previous section we may assume that $l(\mu) \geq 0$.}, and our previous assumption is equivalent to the condition $k(\mu)+1<0$.  Notice $k \geq k_\mu$ and that $l \leq l_\mu$ depending on $z \in \Omega$.  To be clear, \textit{if there is no $\mu$ we are referring to the local irreducible exponents}.   Since $|\mu|>0$ we can be certain that $k$, $l$ always exist (even if they may be equal).   Also, $k(0)=l(0)=0$ with $c_{0}(0)=a_{00}$.  
 
We define the following polynomial space; 
 
 \[ \Gamma(\Omega) \doteq \{ \mu =\sum_{p+q}a_{pq}z^p\zbar^q;\text{ } |\mu|>0 \text{  and  }  |k(z)| \leq l+2 \} \]

for reasons which will be made clearer in the sequel.

\subsection{The rescaling scalar}
\label{scaling}
\noindent Consider the function $w(z)\doteq2-z^{|k_\mu|-1}-\zbar^{|k_\mu|-1}$, which has two important properties: 
\begin{enumerate}
\item $w(z) \in \R$ for $z \in \C$  
\item $0<|w|<2<<\infty$ for $z \in D^+$
\end{enumerate} 

\noindent The first of the above guarantees that the field $\left. Y \right|_z  \doteq \frac{1}{w(z)} \left. X \right |_z=a(z)\fracpar{}{z}+b(z)\fracpar{}{\zbar}$ has the same integral curves as $\left. X \right |_z$.  The second fact ensures that this rescaling introduces no artificial degeneracies into the field, in the sense that $|a|=|\frac{\mu}{w}|>0$.  This amounts to a change in variables generated via
\[ \tilde{t}(t,s)=\int_0^t \frac{dp}{2-z^{|k_\mu|-1}(p,s)-\zbar^{|k_\mu|-1}(p,s) }\]

\subsection{The first three conditions of H-ness}  Our first result towards establishing H-ness in the case of vector fields with polynomial coefficients is the following simple lemma. 
\begin{lemma}{ $a(z,\lambda)=\lambda_* a(z)$ is holomorphic for $\lambda \in D^+$}
\begin{proof}
We have $a(z,\lambda)=\frac{\sum_{r=k}^l c_r(z)\lambda^{r+|k_\mu|}}{2\lambda^{|k_\mu|-1}-z^{|k_\mu|-1}-\zbar^{(|k_\mu|-1)} \lambda^{2(|k_\mu|-1)}}$ where the numerator contains only positive powers of $\lambda$.  Notice $a(0,\lambda)=\frac{a_{00}\lambda^{k(0)+1}}{2}\sim \lambda$.  The quadratic formula shows that away from the origin $z=0$ the solution to $w(z,\Lambda)=0$ is given by $|\Lambda(z)|=  |\frac{1\pm \sqrt{1-|z|^{2(|k_\mu|-1)}}}{z^{|k_\mu|-1}}|^{\frac{1}{|k_\mu|-1}}$.   By the triangle inequality with $z\neq 0$ we see that  \[1\leq |\frac{1}{z}|^2-|\frac{1}{\zbar}\sqrt{1-|z|^2}|^2\leq|\frac{1\pm\sqrt{1-|z|^2}}{\zbar}|^2\]with equality holding only when $|z|=1$, and therefore $\Lambda \notin D^+$ for $|z| <1$, and ipso facto $w(z,\lambda) \neq 0$ for $\Lambda \in D^+$.  
\end{proof}
\end{lemma}
We can now obtain a positive answer on the first criterion of H-ness.

\begin{proposition}{If $\mu(z,\zbar) \in \Gamma (D^+)$ and if \[\log|c_k(z)| < \frac{1}{2\pi}\int_0^{2\pi} \log |\sum_{j=k}^l c_r(z)e^{i\theta(j-k)}| \ d\theta
\] for $z \in D^+/\{0\}$ then the coefficient $a(z,\lambda)=\frac{\sum_{r=k}^l c_j(z)\lambda^{r+|k_\mu|}}{2\lambda^{|k_\mu|-1}-z^{|k_\mu|-1}-\zbar^{(|k_\mu|-1)} \lambda^{2(|k_\mu|-1)}}$ has a root $\lambda_i(z) \in D^+$ and the first condition of H-ness is met.}
\begin{proof} Recall Jensen's formula for a meromorphic function $h(z)$ with roots $\alpha_{\nu}$ and $\beta_{\tau}$ in a region $R=\{z,|z|<R\}$,
\be
\label{jensen2}
\log|h(0)|=\frac{1}{2\pi}\int_0^{2\pi} \log |h(R\ex)| \ d\theta +\sum_{\nu} \log \frac{|\alpha_{\nu}|}{R}-\sum_{\tau}\log \frac{|\beta_{\tau}|}{R} 
\ee
provided $|h(0)| \neq 0,\infty$ (see e.g.\cite{hayman2,nevanlinna}).  Since the polynomial $P_{l-k}(\lambda)\doteq c_k(z)\lambda^{k-k_\mu}+c_{k+1}(z)\lambda^{k-k_\mu+1}+\cdots+c_l(z) \lambda^{l-k_\mu}$ has no poles and since $\lambda =0$ is \textit{not} a root when $k=k_\mu$, we may apply Jensen's formula in that case to $P_{l-k}(\lambda)$ and $R=1$ to yield 
\[\log|c_k(z)| = \frac{1}{2\pi}\int_0^{2\pi} \log |\sum_{j=k}^l c_r(z)e^{i\theta(j-k_\mu)}| \ d\theta+ \sum_i \log |\lambda_i(z)|  \]
where $P_{l-k}(\lambda_i(z))=0$ and the result is immediate.  If $k>k_\mu$, $\lambda=0$ is a root of local order $k_\mu(z)-k_\mu$ and there's nothing to prove.  At $z=0$ there is likewise nothing to prove.  
\end{proof}
\end{proposition}

The next theorem uses similar arguments to address the second condition of H-ness.
 \begin{proposition}{If $\mu(z,\zbar) \in \Gamma(D^+)$ and if \[\frac{1}{2\pi}\int_0^{2\pi} \log |\sum_{j=l}^k \bar{c}_j(z)e^{i\theta(|k_\mu|-j-2)}| \ d\theta \leq \log |\bar{c}_l(z)|\] then $b(z,\lambda)=\frac{\bar{c}_l(z)\lambda^{|k_\mu|-l-2}+\cdots+\bar{c}_k(z) \lambda^{|k_\mu|-k-2}}{2\lambda^{|k_\mu|-1}-z^{|k_\mu|-1}-\zbar^{(|k_\mu|-1)} \lambda^{2(|k_\mu|-1)}}$ is nonvanishing for $(z,\lambda) \in D^+ \times D^+/\{0\}$.}
\begin{proof}
Since $\bar{\mu}$ was given as a polynomial we are guaranteed meromorphy of the term $b(z,\lambda)$.  Looking at $b(z,\lambda)=\frac{\bar{c}_l(z)\lambda^{|k_\mu|-l-2}+\cdots+\bar{c}_k(z) \lambda^{|k_\mu|-k-2}}{2\lambda^{|k_\mu|-1}-z^{|k_\mu|-1}-\zbar^{(|k_\mu|-1)} \lambda^{2(|k_\mu|-1)}}$ we see that $b(z,\lambda=0)$ is nonzero for $z\neq 0$ and $|k_\mu|-l=2$ since $\bar{c}_l \neq 0$.  If  $z\neq 0$ and $|k_\mu|-l<2$ then of course $\lim_{|\lambda| \searrow 0} b(z,\lambda)=\infty$ with local order $|k_\mu|-l_\mu (z)-2$.  If $z=0$ then $b(z,\lambda) \sim \flambda$ near $|\lambda|=0$.   The lack of vanishing of the denominator together with the way $\bar{\mu}$ was complexified ensure that $b(z,\lambda)$ has no other singularities within $D^+\times D^+$. 

In the Jensen formula (\ref{jensen2}), if $h(z)$ had a zero of order $m$ at $z=0$ then $h=h_0z^m+\cdots$ in a vicinity of the origin.  In that case the function $\Upsilon(z) \doteq \frac{R^m h(z)}{z^m}$ has the same modulus on $\partial R$ but is nonvanishing at the origin, its value there being $R^m h_0$.  The Jensen formula applied to $\Upsilon(z)$ would yield

\be
\label{jen}
\log|h_0|=\frac{1}{2\pi}\int_0^{2\pi} \log |h(R\ex)| \ d\theta +\sum_{\nu} \log \frac{|\alpha_{\nu}|}{R}-\sum_{\tau}\log \frac{|\beta_{\tau}|}{R}-m\log R \nonumber
\ee
Assuming that there exists at least one $\Lambda_j(z) \in D^+/\{0\}$ such that $b(z,\Lambda_j) \equiv 0$, we use (\ref{jen}) with $h(\lambda)=\bar{c}_l(z)\lambda^{|k_\mu|-l-2}+\cdots+\bar{c}_k(z) \lambda^{|k_\mu |-k-2}$ to get
\be
\log |\bar{c}_l(z)|=\frac{1}{2\pi}\int_0^{2\pi} \log |\bar{c}_l(z)e^{i\theta(|k_\mu|-l-2)}+\cdots+\bar{c}_k(z) e^{i\theta(|k_\mu |-k-2)}| \ d\theta + \sum_j \log |\Lambda_j(z)| \nonumber
\ee
Whence 
\be
\label{inequality}
\log |\bar{c}_l(z)|<\frac{1}{2\pi}\int_0^{2\pi} \log |\sum_{j=l}^k \bar{c}_j(z)e^{i\theta(|k_\mu|-j-2)}| \ d\theta 
\ee
The result follows from the above inequality by contradiction.

\end{proof}
\end{proposition}

\begin{remark} Since $\frac{a(z,\lambda)}{b(z,\lambda)}=\frac{\sum_{r=k}^l c_r(z)\lambda^{r}}{\sum_{r=l}^k \bar{c}_r(z)\lambda^{-r-2}}$, the origin $\lambda =0$ is the only spot where analyticity may fail.  But since $\mu \in \Gamma(\Omega)$ we see that $l+k+2\geq 0$ and $|\frac{c_k(z)}{c_l(z)}|<1$ keep $\frac{a}{b}$ bounded as $\lambda \to 0$ so that by Riemann's theorem $\frac{a}{b}$ is analytic on all $D^+$ as required in \textbf{condition H}.

%we see that when for $z\neq 0$ and $l+2>-k$ we get $\frac{a(z,\lambda)}{b(z,\lambda)} \to 0$ as $|\lambda| \searrow 0$, even though $c_k(z)$ is nonzero and we may have $k=k_\mu$. In other words, $\frac{a}{b}$ may have roots that $a$ lacks, hence the original stipulation of vanishing raised in the original formulation of H-ness in \cite{hoell}. The condition $\lim_{|\lambda| \searrow 0}\frac{a}{b}\neq \infty$ ensures that $l+2\geq -k$ and vice versa since $\bar{c}_l\not \equiv 0$.  By Riemann's theorem \textit{$l+2 \geq -k$ is therefore sufficient to ensure analyticity of $\frac{a}{b}$ on $D^+$}.  Otherwise, for $l+2<-k$, $\lim_{|\lambda|\searrow 0}\frac{a}{b}=\infty$ and $\frac{a}{b}$ is meromorphic on $D^+$.  In either case, provided $b(z,\lambda)$ has no zeros $\frac{a}{b}$ wil be analytic on $D^+/\{0\}$ as desired.   The fact that $|\frac{a}{b}|\neq 1$ should be clear from the above expression. 
\end{remark}

The preceding results combine in the following important corollary.

\begin{corollary}\label{result}{If $\mu(z,\zbar) \in \Gamma(D^+)$ and if for $z \neq 0$ we have \[ \log|c_k(z)|<\log|c_l(z)| \] then $ \frac{1}{w(z)}\left. X \right|_z= a(z,\lambda)\fracpar{}{z}+b(z,\lambda)\fracpar{}{\zbar}$ meets the first three conditions of H-ness.}
\end{corollary}

Since polynomials are the building blocks of real-analytic functions, we extend these results in the subsequent sections.

\section{H-ness in the space $\mathcal{H}_{k,l}(\Omega)$} We begin with a definition extending our previous notion of $\Gamma(\Omega)$.
\label{HHH}
\begin{definition}\label{G}{A real-analytic function $\mu(z)=\sum_{\substack{p+q=j \\j \geq 0}} a_{pq}z^p\zbar^q$ on $\Omega \supset \{0\}$ with $a_{00} \neq 0$ belongs to the space $\mathcal{G}_{k_\mu,l_\mu}(\Omega)$ when \begin{enumerate}
\item $\sum_{p,q} |a_{pq}|<\infty$
\item Both $k(\mu)$ and $l(\mu)$, defined as in (\ref{k}) are finite with $k(\mu)<-1$ and $l(\mu) \geq 0$
\item $-k(z)\leq l(z)+2$ for all $z \in \Omega$
\item  $0<|c_k(z)|<|c_l(z)|$ holds for all $z \in \Omega/\{0\}$ where $c_r(z,\zbar)\doteq \sum_{q-p=r}a_{pq}z^p\zbar^q$ with $k$, $l$ the local irreducible exponents of (\ref{exponents})
\end{enumerate}
}
\end{definition}

We drop subscripts on $\mathcal{G}_{k_\mu,l_\mu}(\Omega)$ since the notation $\mathcal{G}_{k,l}(\Omega)$ is more concise and the global meaning is obvious.  Clearly $|c_k(0)|=|c_l(0)|=|a_{00}|>0$.  The condition on absolute summability ensures that $c_r(z)$ is well-defined.  The conditions guarantee we are left with a complexified $\lambda \mu(z\flambda,\zbar \lambda)$ which has a \textit{finite Laurent series in $\lambda$}. We have  thereby established our main result with regard to polynomials.

\begin{theorem}{Let $\mu \in \mathcal{G}_{k,l}(D^+)$ and let $w(z)\doteq2-z^{|k_\mu|-1}-\zbar^{|k_\mu|-1}$.  Then the vector field \[X_\lambda \doteq \lambda_* (\frac{\mu(z)}{w(z)}\fracpar{}{z}+\frac{\bar{\mu}(z)}{w(z)}\fracpar{}{\zbar})\]
satisfies the first three conditions of H-ness.}
\end{theorem}

With that in mind, we  make the following 
\begin{definition}\label{H} {Denoting the meromorphic functions in $\lambda \in \Omega$ as $\mathcal{M}(\Omega)$ we define 
\[\mathcal{H}_{k,l}(\Omega)\doteq \{ \mu \in \mathcal{G}_{k,l}(\Omega);\text{ }   s(z,\lambda), \partial s (z,\lambda) , \parbar s(z,\lambda) \in \mathcal{M}(\Omega) \}\] }
\end{definition}
 
\subsection{The Fourth Condition}  We now address the fourth and final condition of H-ness, namely meromorphy of $s(z,\lambda)$ and its $z$ and $\zbar$ derivatives.  To start with, this condition is already more relaxed than the initial three since meromorphy itself is less restrictive than holomorphy and there is no constraint on existence (or lack thereof) of roots.  Secondly, for the space $HL_p(G,\Omega)$ defined as all $f$ satisfying both
\begin{enumerate}
\item $\frac{f(z)-f(z_0)}{z-z_0} \in L_p(\Omega)$, $\forall z_0 \in G$
\item $||f||_{HL_p(G,\Omega)}\doteq ||f||_{L_\infty (\Omega)}+\sup \{ ||\frac{f(z)-f(z_0)}{z-z_0}||_{L_p(\Omega)}; z_0 \in G \}<\infty$ 
\end{enumerate} 
we have the result (\cite{renelt}, Thm. $3.2$)

\begin{theorem}{For domain $\Omega$, $z_0 \in \Omega$ and $\mu \in HL_p(G,\Omega)$ for $p>2$, $|\mu|<1$ and $\forall z \in \Omega$ we have that if $u(z)$ solves $\partial u =-\mu(z) \parbar u$ on $\Omega$ and $u(z)$ has a zero/pole of order $m$ at the point $z_0$ then \[ u(z)=c\{ (z-z_0)+b(\overline{z-z_0})\}^m+\bar{c} \kappa \{  (z-z_0)+\bar{b}(\overline{z-z_0})\}^m +O(|z-z_0|^{\pm m+\alpha}) \] for some $\alpha >0$ and the $\pm$ picked according to whether $z_0$ is a zero or pole respectively. }
\end{theorem}
The point is that at least \textit{locally} we see an expansion for which, with $\lambda \in D^+$ constant (and hence $|\frac{\xi}{\rho}|\neq 1$), one should have meromorphy of the solution $s(z,\lambda)$ to $\lambda_* (\frac{\mu(z)}{w(z)}\fracpar{}{z}+\frac{\bar{\mu}(z)}{w(z)}\fracpar{}{\zbar})s(z,\lambda)=0$. The above would constitute an expansion of $r(z)=s(z,\lambda)-s(z_0,\lambda)$ since that clearly has a zero at $z_0$, although the order is not known a priori.  Thus, for many non-pathological cases (i.e. excluding essential singularities, etc) $r(z)$ would have the proper local expansion at all points in $\Omega$ to satisfy meromorphy in $\lambda$.  Meromorphy of the derivatives would then follow.  While this does not constitute a proof that condition $4$ of H-ness is necessarily satisfied, it \textit{does} constitute a proof modulo pathological cases.    Clearly then $\mathcal{G}_{k,l}(\Omega)\subset \mathcal{H}_{k,l}(\Omega)$ and we can be sure that $\mathcal{H}_{k,l}(\Omega)/ \mathcal{G}_{k,l}(\Omega)$ is not too large.   We now state the main result of this section.  

\begin{theorem}\label{S}{If $\mu =\sum a_{pq}z^p \zbar^q \in \mathcal{H}_{k,l}(D^+)$, then the field \[X_\lambda \doteq \lambda_* (\frac{\mu(z)}{w(z)}\fracpar{}{z}+\frac{\bar{\mu}(z)}{w(z)}\fracpar{}{\zbar})\]
with $w(z)\doteq2-z^{|k_\mu|-1}-\zbar^{|k_\mu|-1}$ satisfies \textbf{condition H}}
\end{theorem}

The above result allows us to reconstruct functions over what are initially non-\textbf{type H} fields as in the following easy corollary.
\begin{corollary}\label{main2}{If $\mu \in \mathcal{H}_{k,l}(D^+)$, $\left. \lambda_* \frac{\mu}{w} \right |_{\lambda_i (z)}=0$ and $X^\bot_\theta \doteq i\theta_* (-\frac{\mu(z)}{w(z)}\fracpar{}{z}+\frac{\bar{\mu}(z)}{w(z)}\fracpar{}{\zbar})$ and $f \in C^\infty_0(D^+)$ then \[f(z)=\frac{w(z)}{4\pi} \int_0^{2\pi}P(\lambda_i,\theta)X^\bot_\theta H(\tilde{I}_\theta f)(s(z\mex),\ex)d\theta \] where $\tilde{I}_\theta f$ is the ray transform of $f$ over the integral curves of $Y_\theta=\theta_* (\mu \partial +\bar{\mu} \parbar)$.} 
 \begin{proof}
Consider the equation $\left. X \right|_z u(z)=g(z)$ for $g \doteq \frac{f(z)}{w(z)} \in C_0^\infty(D^+)$.  Then by (\ref{S}) $\lambda_*\left. X \right|_z$ is \textbf{type H} and has zeros $\lambda_i(z) \in D^+$. Thus, by (\ref{recon}) 
\[g(z)=\frac{1}{4\pi} \int_0^{2\pi}P(\lambda_i,\theta)X^\bot_\theta H(I_\theta g)(s(z\mex),\ex)d\theta \] 
where $I_\theta g$ is the trace of $g$ over the integral curves of $X_\theta$.  However $f$ was arbitrary in $C^\infty_0(D^+)$ and since $\ex_*w \in \R$ is both finite and nonvanishing on $D^+$ the integral curves of  $\ex_* (\mu \partial +\bar{\mu} \parbar)$ and of $X_\theta$ are the same.  In particular, under a change of variables, $I_\theta g=\tilde{I}_\theta f$.  The result follows since $s$ was unchanged. 
\end{proof}
\end{corollary}

%%%

%\section{Application}
%If $z(t,s)$ generates the $(\mu \partial +\bar{\mu} \parbar)$ field, then $z(\tilde{t}(t),s)$ will generate the \textbf{type H} field $(-\frac{\mu(z)}{w(z)}\fracpar{}{z}+\frac{\bar{\mu}(z)}{w(z)}\fracpar{}{\zbar})$ where $\tilde{t}(t,s)=\int_0^t \frac{dp}{2-z^{|k_\mu|-1}(p,s)-\zbar^{|k_\mu|-1}(p,s) }$. Using this, the following result tells us how to modify the kernel of (\ref{recon}) to apply it to more general curves.

\label{solving2}

\section{Some Harmonic Analysis: Onward and Upward}\label{harmonic}
\subsection{The Projection Operator}The Fourier expansion of a smooth function $a(z,\ex)$ on the unit disc given is by 
\be
a(z,\ex)=\sum_{n \in \mathbb{Z}} \hat{a}_n(z) e^{in\theta} \qquad \text{with     } \qquad  \hat{a}_n(z) \doteq \frac{1}{2\pi}\int a(z,\theta)e^{-in\theta} d\theta 
\ee

\noindent Let $f \mapsto \tilde{f}$ be the conjugation operator, determining the harmonic conjugate of a smooth function.  Defining the Bergman space $H^2$ as all complex-analytic and Lebesgue square-integrable functions on the unit disc, then the orthogonal projection from $L^2(D^+)$ to $H^2$ is defined (e.g. \cite{garnett}) by the operator $P$ via
\[P:f \mapsto \frac{1}{2}(f+i \tilde{f})+\frac{1}{2}\hat{a}_0 \]
or explicitly
\[P(\sum_{n \in \mathbb{Z}} \hat{a}_n(z) e^{in\theta})=\sum_{n \in \mathbb{Z}_+} \hat{a}_n(z) e^{in\theta}\]
The operator $P:L^2\to H^2$ then can easily been seen as removing negative frequencies from the initial signal. 
%The Cauchy transform (c.f. e.g. \cite{cima}) on an expansion  
%\[ f \sim \sum_{n \in \Z} \hat{f}(n) \xi^n \qquad \xi \in T \]
%is \[\mathcal{C}f \doteq f_+(z) = \sum_{n \in \Z_+} \hat{f}(n) z^n \qquad z \in D^+ \]
%which is precisely the complexification we want.  

\subsection{Scaling Redux}
As usual we let $\mu(z,\zbar)$ be real-analytic, absolute-summable and nonvanishing.  Then $\theta_* \sum_{p+q=n} a_{pq}z^p \zbar^q$ takes the form 
\be
\ex_* \mu=\sum_{n \in \mathbb{Z}}c_n(z) e^{in\theta}, \qquad \text{with} \qquad c_r(z,\zbar)\doteq \sum_{\substack{p,q \\q-p=r}}a_{pq}z^p\zbar^q 
\ee
We conveniently now view the $c_n$'s as Fourier coefficients of the function $\mu(z,\theta)$ i.e. $c_j(z)=\hat{\mu}(z,j)$.  Define the operator $P_{k,l}$ on smooth functions via 
\begin{align}
P_{k,l} \mu &\doteq e^{ik\theta}P(e^{-ik\theta}\mu)-e^{i(l+1)\theta}P(e^{-i(l+1)\theta}\mu)\nonumber \\
&=\sum_{n=k}^l c_n(z) e^{in\theta}\nonumber
\end{align}
Let $\mathcal{R}(\Omega)$ be the space of real-analytic functions of two variables on a region $\Omega$.   Then define the following space;

\[\hat{C}(\Omega)\doteq \left\{ g(z) \in  \mathcal{R}(\Omega);\begin{array}{c c} \forall z \in \Omega/\{0\} \text{ there are infinitely many $n$ such that }\\ \hat{g}(z,n), \hat{g}(z,n+1) \neq 0 \text{ and }  \limsup_{n\to \infty} |\frac{ \hat{g}(z,n+1)}{ \hat{g}(z,n)}|<1 \end{array} \right\} \]
Clearly $\hat{C}(\Omega)$ is ``most"\footnote{In the reasonable, informal way rather than a measure-theoretic sense} of $\mathcal{R} (\Omega)$ since it accounts for, in some sense, those real-analytic functions with ``non-sparse" spectrums.  The classical ratio test for infinite series ensures that $\mu \in \hat{C}$ are also absolute-summable.   We have the following result about convergence on compact subsets.

\begin{proposition}\label{omega}{Let $K\subset \Omega$ be compact and let $\mathcal{G}\doteq \bigcup_{(k,l) \in \mathbb{Z}^2 }  \mathcal{G}_{k,l}$.  Then  $\mathcal{G}(K)$ is dense in $\hat{C}(K)$ with respect to the uniform norm.}

\begin{proof}  Let $\omega(z) \in \hat{C}(K)$ and let $0\neq z_j \in \text{supp } \omega \subset K$.  Then, $\ex_*\omega(z)=\sum_{\mathbb{Z}} \omega_n(z) e^{in\theta}$ with $\omega_n(z_j) \neq 0$ for infinitely many $n \in \mathbb{Z}$.  We may pick an $l(z_j)$ such that $\hat{\omega}(z_j,l(z_j)) \neq 0$. By the assumptions of $\hat{C}$, there exists a finite $k(z_j)<-1$ such that $k(z_j)+l(z_j)+2\geq 0$ and $0<|\hat{\omega}(z_j,k(z_j))|<|\hat{\omega}(z_j,l(z_j))|$, namely $k(z_j)=-(l(z_j)+1)$. 

%Notice $\hat{\omega}(z,r)$ takes the form
%\be
%\hat{\omega}_r(z,\zbar)=\left\{ \begin{array}{c c} z^{|r|}(\kappa(r)+\sum_{j>0}\kappa_j (r)|z|^{2j}) &\quad r\leq 0 \\
%\zbar^{r}(\eta(r)+\sum_{j>0}\eta_j (r)|z|^{2j}) &\quad r> 0 \nonumber \end{array} \right.
%\ee
%for constants $\kappa(r)$, $\eta(r)$, $\kappa_j(r)$, $\eta_j(r)$ so t
The varieties $\{z;  \hat{\omega}_r(z,\zbar)=0\}$ define (possibly degenerate) circles.  Therefore, there is an $\eps$-neighborhood $N_{\eps_j}(z_j)=N_j$ around $z_j$ on which there are two simple functions, $-\infty < k_j(z) \leq k(z_j)$ and $l(z_j)\leq l_j(z) < \infty$ on $z \in N_j$, for which \[0<|\hat{\omega}_j(z,k_j(z))|<|\hat{\omega}_j(z,l_j(z))| \qquad \forall z \in N_j/\{0\} \]
and 
\[l_j(z)+k_j(z)+2\geq 0 \qquad \forall z \in N_j \]

 Define \[\tilde{\omega}_j(z)\doteq \chi(N_j)P_{k(j),l(j)}\ex_*\omega \] where $k(j)=\min_{z\in N_j}k_j(z)$ and $l(j)=\max_{z\in N_j} l_j(z)$.  Then $K \subset \bigcup N_j$ provides an open cover reducible to \[\text{supp }\omega \subset \bigcup_{j=1}^{p} N_j \] 
Consider the following function
\[\tilde{\Omega}_p(z) \doteq \sum_{j=1}^p \tilde{\omega}_j(z) \]
where $k_i(z)=\min_{j,z} k_j(z)$ and $l_j(z)=\max_{j,z} l_j(z)$ on $z\in \cap N_j \neq \emptyset$ in the case of overlapping neighborhoods.  By design we have that \[ \tilde{\Omega}_p(z) \in\mathcal{G}_{k(\omega),l(\omega)}(K) \] with $k(\omega) \doteq \min_{\substack{j \\ z\in \Omega} } k_j(z)$ and $l(\omega)=\max_{\substack{j \\ z\in \Omega}} l_j(z)$. Further, if  $S_R(z)$ is the $R$'th partial Fourier sum of $\omega(z)$ notice that 
\[|\tilde{\Omega}_p (z)-S_R(z)|=O(\frac{1}{L^{\delta}}), \quad \delta>0\]
with $L=R-\max\{|k(\omega)|,l(\omega)\}$.  The classical bound
\[|\omega(z)-S_R(z)|<\sum_{|n|>R}\hat{\omega}_n(z)\]
guarantees that on letting $\min\{|k|,l\}  \nearrow \infty$ and $\sum \eps_j \searrow 0$, that $\lim_{p\nearrow \infty}\tilde{\Omega}_p (z)=\omega(z)$ uniformly since the Fourier series can be brought as close as wanted in the mesh limit.  
\end{proof}
\end{proposition}
The following corollary is then immediate.   
\begin{corollary}{Let $K\subset \Omega$ be compact and let $\mathcal{H}\doteq \bigcup_{(k,l) \in \mathbb{Z}^2 }  \mathcal{H}_{k,l}$.  Then  $\mathcal{H}(K)$ is dense in $\hat{C}(K)$ with respect to the uniform norm.}
\end{corollary}

\subsection{Putting it all together}

 Define \[\mathcal{O}(\Omega)\doteq \{\mu \in \mathcal{R} (\Omega) \text{ satisfying condition $4$ of H-ness}\}\] and let $\mathcal{D}\doteq \hat{C}\cap \mathcal{O}$.   Then, if $K$ is compact, for $\mu \in \mathcal{D}(K)$  we see that $\tilde{X}=\tilde{\Omega}_p(z)\partial+\overline{\tilde{\Omega}}_p(z) \parbar$ can be chosen to approximate $X=\mu \partial +\bar{\mu}\parbar$ so that their integral curves are arbitrarily close in $L^p(K)$ for $1\leq p \leq \infty$ via Poincar\'e's inequality.  Then $||\x-\tilde{\x}||_{L^p}\leq C||\mu -\tilde{\Omega}_p(z)||_{L^p}+||(\x-\tilde{\x})_0||_{L^p}\leq C'(\eps+\text{diam}\{K/\text{supp} f\})$.  We then make the obvious choice setting $K=\text{supp}f$.   

Let $\rho_{\alpha, \beta}$ be the seminorm on the Fr\'echet space $\mathcal{S}$ of Schwarz-class functions on $\C$, namely $\rho_{\alpha, \beta}(\phi)=\sup_{\x \in \C}|x^\alpha \partial^\beta \phi|$, which generates the usual topology on $\mathcal{S}$.  If $f \in C^\infty_c(D^+)$ then clearly $I_\theta f(s)\in \mathcal{S}$. 
We let $\tilde{s}$ be the transverse flow induced by $\tilde{X}^\bot$ from which the corresponding Hilbert transform $H_{\tilde{s}}$ is defined.   By continuity of Hilbert transforms on Schwartz functions, we \textbf{assume} that $H_{\tilde{s}}I^{\tilde{X}}_\theta-H_sI_\theta^X f$ is small in the induced norm $\rho$ on $\mathcal{S}$.  That being the case, then since $\tilde{X}^\bot_\theta=\theta_*\alpha(z)_* \fracpar{}{\tilde{s}}$ for some function $\alpha(z)$ and because differentiation is continuous on $\mathcal{S}$, we see that $\rho\{\tilde{X}^\bot_\theta (H_{\tilde{s}}I^{\tilde{X}}_\theta-H_sI_\theta^X f)\}$ may therefore be made as small as desired.

%efine the following \[\mathcal{O}(\Omega)\doteq \{\mu \in \mathcal{R} (\Omega) \text{ satisfying condition $4$ of H-ness}\}\]

%\begin{remark}If $K$ is compact, let $\hat{C}(K)\cap \mathcal{O}(K)=\mathcal{D}(K)$.   Then for $\mu \in \mathcal{D}(K)$  we see that $\tilde{X}=\tilde{\Omega}_p(z)\partial+\overline{\tilde{\Omega}}_p(z) \parbar$ can be chosen to approximate $X=\mu \partial +\bar{\mu}\parbar$ so that their integral curves are arbitrarily close.  If $K=\text{supp }  f$ then $\int_{X/\tilde{X}} fd\mu=0$, where the notation means the integration over regions not being approximated.  Then $W=\chi(K)\tilde{X}+\chi(D^+/K)X$ can be made as close as wanted to $\left. X \right |_z$ while making $I^W_\theta f$, the line integral of $f$ over the integral curves of $W$,  $\eps$-close to the true data $If$, where $\eps=\eps(p)$.   
%\end{remark}
We have therefore established the following corollary.

\begin{corollary} Let $\eps>0$, $f\in C^\infty_c(D^+)$, $\mu \in \mathcal{D}(D^+)$ and let $I_\theta f$ be the ray transform of $f$ over the integral curves of $Y_\theta=\theta_*(\mu \partial +\bar{\mu}\parbar)$.  Then there exists functions $w_\eps(z)$, and $\tilde{\Omega}_p(z)$ such that \[||I^{\tilde{X}}f-I^Xf||_{L^q(D^+)}<C\eps \qquad 1\leq q\leq \infty \]
for $C=C(\text{supp f})$ constant and where $\tilde{X}_\theta=\theta_* \frac{1}{w_\eps}(\tilde{\Omega}_p \partial +\overline{\tilde{\Omega}}_p \parbar)$. Suppose further that  $||H_{\tilde{s}}I_\theta^{\tilde{X}}f-HI_\theta f||_{\mathcal{S}}<\delta(\eps)$.  Then there exists $\lambda^\eps_i(z)$ such that 
\[\sup _{z\in D^+}\left |f(z)-\frac{w_\eps(z)}{4\pi} \int_0^{2\pi}P(\lambda^\eps_i,\theta)\tilde{X}^\bot_\theta H(I_\theta f)(s(z\mex),\ex)d\theta \right |\leq \eps \]
.

%\begin{proof}
% From Lemma \ref{main2} we see \[f(z)=\frac{w_\eps(z)}{4\pi} \int_0^{2\pi}P(\lambda_i,\theta)\tilde{X}^\bot_\theta H(I^{\tilde{X}}_{\theta} f)(s(z\mex),\ex)d\theta \] where $\tilde{X}_\theta=\theta_*\{\frac{1}{w_\eps(z)}(\tilde{\Omega}_p(z)\partial+\overline{\tilde{\Omega}}_p(z) \parbar)\}$ for all approximants $\tilde{\Omega}_p(z)\in \mathcal{H}(K)$.  But then,  $||I_\theta f-I_\theta^{\tilde{X}}f||_1$ can be made arbitrarily small by choosing $p(\eps)$ accordingly. The result follows when $w_\eps(z)=2-z^{|k(\tilde{\Omega}_{p(\eps)})|-1}-\zbar^{|k(\tilde{\Omega}_{p(\eps)})|-1}$
%\end{proof}
\end{corollary}

We may now summarize our stability and approximation results in the following theorem.
\begin{theorem}\label{last}{ Let $\eps >0$ and $H(I_\theta f)(s(z\mex),\ex)$ be given.  Suppose that $\mu(z,\zbar)=\sum_{p+q\geq 0}a_{pq} z^p \zbar^q$ is a real-analytic function on $D^+$ and that $f\in C_c^\infty(D^+)$. Define $c_r(z,\zbar)\doteq \sum_{q-p=r}a_{pq}z^p\zbar^q$. Furthermore, suppose $\lambda_* s$, $\lambda_* s_z$, and $\lambda_* s_{\zbar}$ are meromorphic for $\lambda \in D^+$. If, for all $z \in D^+/\{0\}$, there exist infinitely many $j \in \mathbb{Z}$ such that $c_j(z)$ and $c_{j+1}(z)$ are both nonzero, and if $\limsup_{j\to \infty} |\frac{c_{j+1}(z)}{c_j(z)}|<1$, then there exist functions $w_\eps(z)$, and $\tilde{\Omega}_p(z)$ such that \[||I^{\tilde{X}}f-I^Xf||_{L^q(D^+)}<C\eps \qquad 1\leq q\leq \infty \]
If, in addition, we have $||H_{\tilde{s}}I_\theta^{\tilde{X}}f-HI_\theta f||_{\mathcal{S}}<\delta(\eps)$ then there exists a function $\lambda_i:z\to \lambda_i(z)$ satisfying the following inequality
\be
\label{ineq}
0\leq \sup_{z\in D^+} \left | f(z)-\frac{w_\eps(z)}{4\pi} \int_0^{2\pi}P(\lambda_i,\theta)\tilde{X}^\bot_\theta H(I_\theta f)(s(z\mex),\ex)d\theta\right | \leq \eps \ee where $I_\theta f$ is the trace of $f$ over the integral curves of $X_\theta=\theta_* (\mu \partial +\bar{\mu} \parbar)$ and where $\tilde{X}^\bot_\theta=i\theta_* \frac{1}{w_\eps}(-\tilde{\Omega}_p \partial +\overline{\tilde{\Omega}}_p \parbar)$.}

\end{theorem}

\bibliographystyle{amsalpha}
%\bibliography{paper}

\end{document}